\documentclass{article}

\PassOptionsToPackage{bookmarks=true}{hyperref}
\usepackage[T1]{fontenc}
\usepackage{layout}
\usepackage[left=2cm,right=2cm,top=2cm,bottom=2cm]{geometry}
\usepackage[geometry]{ifsym}
\usepackage{mathtools}
\usepackage{wrapfig}
\usepackage{centernot}
\usepackage{enumitem}
\usepackage{etoolbox}
\usepackage{amsmath}
\usepackage{svg}
\usepackage{amsfonts}
\usepackage{amsthm}
\usepackage{amssymb}
\usepackage{booktabs}
\usepackage[numbers,sort&compress]{natbib}
\usepackage{bookmark}
\usepackage{float}
\usepackage{titlesec}
\usepackage{tikz-cd}
\usepackage{tabularx}
\usepackage{mathrsfs}
\usepackage{bbm}
\usepackage{esvect}
\usepackage{bbold}
\usepackage{suetterl}
\usepackage{fancyhdr}
\usepackage{helvet}
\usepackage{braket}
\usepackage{graphicx}
\usepackage{subcaption}
\usepackage{tikz, tikzscale, tikz-network}
\usetikzlibrary{external}
\usepackage{standalone}
\usepackage{textcomp}

\usepackage{subfiles}

\newcommand{\beginProofOfTheoremClickable}[1]{\begin{proof}[Proof of Theorem \ref{#1}]}
\newcommand{\beginProofOfLemmaClickable}[1]{\begin{proof}[Proof of Lemma \ref{#1}]}
\newcommand{\beginProofOfCorollaryClickable}[1]{\begin{proof}[Proof of Corollary \ref{#1}]}
\newcommand{\beginProofOfClaimClickable}[1]{\begin{proof}[Proof of Claim \ref{#1}]}
\newcommand{\beginProofOfObservationClickable}[1]{\begin{proof}[Proof of Observation \ref{#1}]}
\newcommand{\beginProofOfExampleClickable}[1]{\begin{proof}[Proof of Example \ref{#1}]}

\newcommand{\beginProofOfTheorem}[1]{\begin{proof}[Proof of Theorem \ref*{#1}]}
\newcommand{\beginProofOfLemma}[1]{\begin{proof}[Proof of Lemma \ref*{#1}]}
\newcommand{\beginProofOfCorollary}[1]{\begin{proof}[Proof of Corollary \ref*{#1}]}
\newcommand{\beginProofOfClaim}[1]{\begin{proof}[Proof of Claim \ref*{#1}]}
\newcommand{\beginProofOfObservation}[1]{\begin{proof}[Proof of Observation \ref*{#1}]}
\newcommand{\beginProofOfExample}[1]{\begin{proof}[Proof of Example \ref*{#1}]}

\newcommand{\NN}{\mathbb{N}}

\newcommand{\ZZ}{\mathbb{Z}}

\theoremstyle{definition}
\newtheorem{Theorem}{Theorem}

\theoremstyle{definition}
\newtheorem{Observation}{Observation}

\theoremstyle{definition}

\theoremstyle{definition}

\theoremstyle{definition}
\newtheorem{Definition}{Definition}

\theoremstyle{definition}
\newtheorem{Claim}{Claim}

\theoremstyle{definition}
\newtheorem{Lemma}{Lemma}

\theoremstyle{definition}
\newtheorem{Corollary}{Corollary}

\theoremstyle{definition}

\theoremstyle{definition}

\newenvironment{proofref}[2]{%
  \begin{proof}[Proof of #1~\ref{#2}]%
}{%
  \end{proof}%
}

\title{Bipartite Turán problem on cographs}
\author{Jakob Paul Zimmermann}
\date{December 2025}

\begin{document}

\maketitle

\begin{abstract}
A \emph{cograph} is a graph that contains no induced path $P_4$ on four vertices or equivalently a graph that can be constructed from vertices by \emph{sum} and \emph{product} operations.
We study the bipartite Tur\'an problem restricted to cographs, also known as the \emph{Zarankiewicz problem on cographs}: for fixed integers $s \leq t$, what is the maximum number of edges in an $n$-vertex cograph that does not contain $K_{s,t}$ as a subgraph?
This problem falls within the framework of induced Tur\'an numbers $\text{ex}(n, \{F, H\text{-ind}\})$ introduced by Loh, Tait, Timmons, and Zhou.

Our main result is a \emph{Pumping Theorem}: for every $s\le t$ there exists a period $R$ and \emph{core} cographs such that for all sufficiently large $n$ an extremal cograph is obtained by repeatedly \emph{pumping} one designated pumping component inside the appropriate core (depending on $n\bmod R$). We determine the linear coefficient of $\text{ex}(n, \{K_{s,t}, P_4\text{-ind}\})$ to be $s-1 + \frac{t-1}{2}$, thereby solving the Zarankiewicz problem on cographs up to a constant. Moreover, the pumping components are $(t-1)$-regular and have $s-1$ common neighbours in the respective core graphs, giving the extremal cographs a particularly rigid extremal star-like shape.

Motivated by the rarity of complete classification of extremal configurations, we completely classify all $K_{3,3}$-free extremal cographs by proof. We also develop a dynamic programming algorithm for enumerating extremal cographs for small $n$.
\end{abstract}

\section{Introduction}

\subsection{Background}

The Tur\'an problem is one of the central questions in extremal graph theory: given a graph $H$, what is the maximum number of edges $\text{ex}(n, H)$ in an $n$-vertex graph that does not contain $H$ as a subgraph?
Tur\'an's Theorem~\cite{Shapira2016} determines this exactly when $H = K_{r+1}$ is a clique: $\text{ex}(n, K_{r+1}) \leq \left(1 - \frac{1}{r}\right)\frac{n^2}{2}$, with the unique extremal graph being the \emph{Tur\'an graph} $T_r(n)$---the complete $r$-partite graph with parts as equal as possible.
Notably, $T_r(n)$ as a complete multipartite graph is a cograph, so in the case of cliques the extremal construction already lies within the class of cographs.
For general graphs $H$, the Erd\H{o}s-Stone-Simonovits Theorem~\cite{Shapira2016} provides the asymptotic answer: $\text{ex}(n, H) = (1 - 1/(\chi(H)-1) + o(1)) \binom{n}{2}$.
However, when $H$ is bipartite (i.e., $\chi(H) = 2$), this Theorem only yields $\text{ex}(n, H) = o(n^2)$, leaving the exact order of magnitude as a challenging open problem~\cite{Furedi2013}.

For complete bipartite graphs $K_{s,t}$ with $s \leq t$, the Turán problem is known as the \emph{Zarankiewicz problem}, named after the Polish mathematician who posed it in the 1950s~\cite{Smorodinsky2023}.
The fundamental upper bound was established by K\H{o}v\'ari, S\'os, and Tur\'an~\cite{KovariSosTuran1954}:
\begin{align}
    \text{ex}(n, K_{s,t}) \leq \tfrac{1}{2}(t-1)^{1/s} n^{2-1/s} + \tfrac{1}{2}(s-1)n. \label{eq:KST_background}
\end{align}

Loh, Tait, and Timmons~\cite{Loh2016} introduced the \emph{induced Tur\'an number} $\text{ex}(n, \{H, F\text{-ind}\})$, defined as the maximum number of edges in an $n$-vertex graph that is $H$-free and induced $F$-free.
Our problem falls within this framework: we study $\text{ex}(n, \{K_{s,t}, P_4\text{-ind}\})$, the bipartite Tur\'an problem restricted to cographs.

\subsection{Cographs}

A \emph{cograph} (short for \emph{complement-reducible graph}) is a graph that contains no induced path on four vertices~\cite{Corneil1981}.
Equivalently, cographs are precisely the graphs that can be constructed from single vertices using two recursive operations: the \emph{product} (also called \emph{join}) and the \emph{sum} (disjoint union).
Both operations are associative and commutative.

For two graphs $G_1, G_2$, define their \emph{product} $G_1 \times G_2$ as their disjoint union together with all edges between $V(G_1)$ and $V(G_2)$:
\[
E(G_1 \times G_2) = E(G_1) \cup E(G_2) \cup \{(v_1, v_2) : v_1 \in V(G_1), v_2 \in V(G_2)\}.
\]
Define their \emph{sum} $G_1 + G_2$ as simply their disjoint union.
Iterating these operations for graphs $G_1, \ldots, G_n$, we define $\prod_{j \in [n]} G_j$ and $\sum_{j \in [n]} G_j$ respectively.

\subsection{Notation}

We collect here the notation used throughout the paper.

\paragraph{Sets and intervals.}
For $n \in \NN$, we write $[n] \coloneqq \{1, 2, \ldots, n\}$ for the set of the first $n$ positive integers.

\paragraph{Graphs.}
For a graph $G$, we denote its vertex set by $V(G)$ and its edge set by $E(G)$.
The \emph{vertex count} is $|G| \coloneqq |V(G)|$ and the \emph{edge count} is $\|G\| \coloneqq |E(G)|$.
For two disjoint vertex sets $A, B \subseteq V(G)$, we write $\|A, B\|$ for the number of edges between $A$ and $B$, i.e.,
\[
\|A, B\| \coloneqq \left| \{ \{u, v\} \in E(G) : u \in A, v \in B \} \right|.
\]
The \emph{neighborhood} of a vertex $v$ in $G$ is denoted $\mathscr{N}(v)$ or $\mathscr{N}_G(v)$ when context requires.
For a vertex set $W \subseteq V(G)$, we write $\mathscr{N}_W(v) \coloneqq \mathscr{N}(v) \cap W$ for the neighborhood of $v$ restricted to $W$. We denote the graph complement of $G$ by $\overline{G}$.

\paragraph{Cliques, bicliques, and edgeless graphs.}
For $r \in \NN$, the \emph{clique} $K_r$ is the complete graph on $r$ vertices.
For $s, t \in \NN$, the \emph{biclique} (or \emph{complete bipartite graph}) $K_{s,t}$ is the graph with vertex partition into sets of sizes $s$ and $t$, where every vertex in one part is adjacent to every vertex in the other part. In case $s = 1$ we call $K_{1,t}$ a \emph{star}.
We denote by $E_r$ the \emph{empty graph} (or \emph{independent set}) on $r$ vertices, i.e., the graph with $r$ vertices and no edges.
The \emph{null graph} $\mathbf{0}$ is the unique graph with no vertices (and hence no edges).
The null graph serves as the identity element for both sum and product: $G + \mathbf{0} = G$ and $G \times \mathbf{0} = G$ for any graph $G$.

\paragraph{Subgraphs and forbidden patterns.}
A graph $G$ contains a graph $H$ as a \emph{subgraph} if there exists an injective map $\varphi: V(H) \to V(G)$ such that $\{u, v\} \in E(H)$ implies $\{\varphi(u), \varphi(v)\} \in E(G)$.
We say $G$ contains $H$ as an \emph{induced subgraph} if additionally $\{\varphi(u), \varphi(v)\} \in E(G)$ implies $\{u, v\} \in E(H)$, i.e., $\varphi$ preserves both edges and non-edges.
A graph $G$ is \emph{$H$-free} if it does not contain $H$ as a subgraph.
A graph $G$ is \emph{induced $H$-free} if it does not contain $H$ as an induced subgraph.

\paragraph{Tur\'an numbers.}
The \emph{Tur\'an number} $\text{ex}(n, H)$ is the maximum number of edges in an $n$-vertex $H$-free graph.
The \emph{induced Tur\'an number} $\text{ex}(n, \{F, H\text{-ind}\})$ is the maximum number of edges in an $n$-vertex graph that is $F$-free and induced $H$-free.
More generally, for a host graph $G$, we define $\text{ex}(G, H)$ as the maximum number of edges in an $H$-free subgraph of $G$, and $\text{ex}(G, \{F, H\text{-ind}\})$ as the maximum number of edges in a subgraph of $G$ that is $F$-free and induced $H$-free.

\paragraph{Edge contribution.}
For a vertex set $X \subseteq V(G)$, we define its \emph{edge contribution} by
\[
e(X) \coloneqq \left\| G[X] \right\| + \left\| X, V(G) \setminus X \right\|,
\]
where $G[X]$ denotes the induced subgraph on $X$.
Intuitively, $e(X)$ counts both the internal edges within $X$ and the cross edges connecting $X$ to the rest of the graph.

\subsection{Main results} \label{sec:main_results_intro}

For $s, t \in \NN$ with $s \leq t$, we call a cograph \emph{$(s,t)$-extremal} if it is edge-maximal among cographs of the same vertex count not containing $K_{s,t}$.
To built some first intuition about the problem, we developed a dynamic programming framework (Section~\ref{subsec:dynamic_programming}) to enumerate extremal cographs for small $n$. Check out the presentation of found extremal $(s,t)$-extremal cographs at \url{https://extremal-cographs.fly.dev}. All code is available at \url{https://github.com/JayPiZimmermann/zarankievicz_cographs.git}
\begin{Definition}[Pumping]\label{def:pumping}
Let $G$ be a graph and let $H \subseteq G$ be a subgraph. We say that we $k$ times \emph{pump} (or $k$ times \emph{pump up}) $H$ in $G$ to obtain the graph $G'$ constructed by adding $k$ copies of $H$, where each copy has every vertex connected to the $V(G - H)$-neighborhood of its corresponding vertex in $H$. Formally:
\begin{align}
    G' \coloneqq \left( V(G) \setminus V(H) \cup \bigcup_{j \in [k+1]} V(H_j), \ E(G - H) \cup \bigcup_{j \in [k+1]} \left( E(H_j) \cup \bigcup_{v \in V(H)} \{ v' \} \times \mathscr{N}_{V(G-H)}(v) \right) \right),
\end{align}
where $H_1, \dots, H_{k+1}$ are copies of $H$ and vertex $v$ in $V(H)$ corresponds to vertex $v'$ in $V(H_j)$.
\end{Definition}

Our first main result shows that for large $n$, extremal cographs arise by pumping up a specific component inside a core graph.

\begin{Theorem}[Pumping Theorem] \label{thm:pumping}
    Let $s, \ t \in \NN$ with $2 \leq s \leq t$.
    There exist $N_{s,t}, R \in \NN$ and core graphs $G_0, \ldots, G_{R-1}$ each of size at most $N_{s,t}$ with designated $(t-1)$-regular components $H_j \subseteq G_j$ with $s-1$ common neighbors outside of $H_j$ such that for any $n$ large enough there is an $(s,t)$-extremal cograph on $n$ vertices obtained by pumping $H_{(n \bmod t)}$ inside of $G_{(n \bmod t)}$.
    Moreover, there are constants $a_0, \ldots, a_{R-1} < 0$ such that for large $n \in \NN$
    \begin{align}
        \text{ex}(n, \{K_{s,t}, P_4\text{-ind}\}) = a_{(n \bmod R)} + \left( \frac{t-1}{2} + s-1 \right) \cdot n.
    \end{align}
    There exists at least one residue class $g \in \mathbb{Z} / R \mathbb{Z}$ such that the $(s,t)$-extremal cographs obtained by pumping $H_q$ in $G_q$ are connected.
    Moreover, for any $n \in \NN$
    \begin{align}
        \text{ex}(n, \{K_{s,t}, P_4\text{-ind}\}) < \left( \frac{t-1}{2} + s-1 \right) \cdot n.
    \end{align}
\end{Theorem}

The key structural tool underlying the Pumping Theorem is the \emph{Lifting Theorem} (Theorem~\ref{thm:lifting}), which shows that extremal cographs admit a specific decomposition into products of small cores with sums of smaller cographs. This lifting procedure controls the connectivity components of extremal cographs and constrains their structure, enabling the eventual pumping behavior.
In Theorem \ref{thm:pumping_profile} we also prove a pumping-like behaviour of cographs edge-maximal with respect to not containing any of a set of bicliques, which is the basis for proving Theorem \ref{thm:pumping}.

\begin{Corollary} \label{thm:balanced_bicliques}
For any $n \in \NN$ and $t \coloneq \lfloor \frac{n}{6} + 1 \rfloor $ and any cograph $G$ on $n$ vertices, either $K_{t,t} \subseteq G$ or $K_{t,t} \subseteq \overline{G}$.
\end{Corollary}

\beginProofOfCorollaryClickable{thm:balanced_bicliques}

By Theorem \ref{thm:pumping} we have $\text{ex}(n, \{K_{s,t}, P_4\text{-ind}\}) < \frac{3(t-1)}{2} \cdot n \leq \frac{n(n-1)}{4}$.
However, by the pigeon hole principle either $G$ or $\overline{G}$ has at least $\frac{n(n-1)}{4}$ edges.
\end{proof}

We remark that this shows exactly that cographs have the \emph{Strong Erd\H{o}s-Hajnal propery} \citep{Cho2024}. However, they proved the existence of a balanced biclique with partite sets of size $\frac{1}{4}n$, yielding a stronger factor as an application of the pigeon hole principle together with upper bounds for the extremal function.

In the following subcases, we can fully classify all extremal cographs - also for small $n$. For forbidden stars $K_{1,t}$, the extremal cographs are regular in most cases.

\begin{Theorem}[Forbidden stars] \label{thm:pumping_stars}
    Let $t,n \in \NN$ with $n \geq t$. In case $n$ is even or $t-1$ is even and $2t \neq n+1$, any $(1,t)$-extremal cograph on $n$ vertices is $(t-1)$-regular.
    Otherwise, it is the sum of a $(t-1)$-regular cograph and a connected remainder of size at most $2t - 3$.
\end{Theorem}

For $K_{2,t}$, the extremal cographs have a particularly simple structure:

\begin{Theorem} \label{thm:classification_2_t}
    Let $t \in \{2,3\}$ and let $G$ be a $(2,t)$-extremal cograph on $n \geq 2$ vertices.
    Then, $G$ is the product of a vertex with another cograph $G'$ that is $(1,t)$-extremal.
\end{Theorem}

However, we remark that for larger $t$ there are $n \leq 2(t-1)$ such that the $(2,t)$-extremal cographs on $n$ vertices do not have a complete vertex. For example, $G \coloneq 3 \cdot K_2 \times 3 \cdot K_2$ is $(2,7)$-extremal on $12$ vertices. However, $\| G \| = 42 < 48 = \left( 2-1 + \frac{7-1}{2} \right)$ so this could not be a pumping component, aligning with the Pumping Theorem.

Our final main result provides a complete classification in the fundamental case $K_{3,3}$:

\begin{Theorem}\label{thm:classification_3_3}
    Let $G$ be an edge-maximal cograph on $n \geq 2$ vertices that does not contain $K_{3,3}$. Then, $G$ is the product of a two-vertex cograph $G_1$ with another cograph $G_2$. Moreover, there is always an edge-maximal cograph on $n \geq 2$ vertices such that $G_1$ is an edge and $G_2$ is the sum of cliques.
\end{Theorem}

The proofs of these Theorems are presented in Sections~\ref{sec:biclique_profiles} and~\ref{sec:main_results}.

\subsection{Further background and related work}

\paragraph{More insights about the Zarankievicz problem.}
The K\H{o}v\'ari-S\'os-Tur\'an bound (\ref{eq:KST_background}), known as the \emph{KST bound}, is known to be tight up to constant factors in several cases.
The cases $s = t = 2$ and $s = t = 3$ were resolved by classical constructions~\cite{Furedi2013}.
Koll\'ar, R\'onyai, and Szab\'o~\cite{KollarRonyaiSzabo1996} introduced \emph{norm graphs} to show tightness when $t \geq s! + 1$, which Alon, R\'onyai, and Szab\'o~\cite{AlonRonyaiSzabo1999} improved to $t \geq (s-1)! + 1$.
Conlon~\cite{Conlon2022_someRemarksZarankievicz} developed a quantitative variant of the random algebraic method to prove that $\text{ex}(K_{n,m}, K_{s,t}) = \Omega(mn^{1-1/s})$ for any $m \leq n^{t^{1/(s-1)}/s(s-1)}$, showing tightness of the KST bound over a broader range than previously known.
However, determining the exact order of magnitude for general $s$ and $t$ remains a central open problem~\cite{Furedi2013, Smorodinsky2023}.
For small parameters, Tan~\cite{Tan2022_zarankieviczSAT} used SAT solvers to compute exact values of the Zarankiewicz function, correcting errors in earlier hand-computed tables and extending the known range of parameters.

\paragraph{Structural constraints and improved bounds.}
A key insight in modern extremal combinatorics is that additional structural constraints on the host graph can yield dramatically improved bounds.
A foundational result in this direction is due to Bonamy, Bousquet, Pilipczuk, Rz\k{a}\.{z}ewski, Thomass\'e, and Walczak~\cite{Bonamy2020}, who proved that $K_{\ell,\ell}$-free graphs avoiding induced $P_5$ have degeneracy $O(\ell^3)$.
Since cographs (graphs avoiding induced $P_4$) are a subclass of induced $P_5$-free graphs, this immediately implies that $K_{\ell,\ell}$-free cographs have at most $O(\ell^3 n)$ edges---a linear bound in $n$.
More generally, Bonamy et al.\ showed that $K_{\ell,\ell}$-free graphs avoiding long induced paths or cycles have degeneracy polynomial in $\ell$, with specific exponents depending on the forbidden pattern.
In a major breakthrough, Nguyen~\cite{Nguyen2025} resolved a 1985 problem of Gy\'arf\'as by proving that induced $P_5$-free graphs are polynomially $\chi$-bounded: there exists $d \geq 2$ such that every induced $P_5$-free graph $G$ satisfies $\chi(G) \leq \omega(G)^d$.

F\"uredi~\cite{Furedi1991} and Sudakov--Tomon~\cite{Sudakov2020} developed powerful techniques for bipartite Tur\'an problems using \emph{neighborhood intersection arguments}, where vertex sets in extremal structures cannot share too many common neighbours.
F\"uredi~\cite{Furedi1991} used this to show that graphs avoiding $L^k$ (the lowest three levels of the Boolean lattice) have $O(n^{3/2})$ edges, confirming a conjecture of Erd\H{o}s.
Sudakov and Tomon~\cite{Sudakov2020} extended these ideas using hypergraph methods: given a dense bipartite graph, they construct a $t$-uniform hypergraph on common neighborhoods and apply the Hypergraph Removal Lemma to embed forbidden subgraphs.

The survey by Keller and Smorodinsky~\cite{Smorodinsky2023} develops a unified approach to the Zarankiewicz problem via \emph{$\varepsilon$-$t$-nets}---a generalization of classical $\varepsilon$-nets where one seeks a small family of $t$-tuples hitting all large hyperedges.
Their key result shows that $K_{t,t}$-free bipartite graphs with VC-dimension $d$ have $O(n^{2-1/d})$ edges, recovering the Fox-Pach-Sheffer-Suk-Zahl bound~\cite{Smorodinsky2023} with a simpler proof.
For geometric intersection graphs, they obtain much stronger bounds: intersection graphs of pseudo-discs have only $O(t^6 n)$ edges when $K_{t,t}$-free (a \emph{linear} bound), and axis-parallel rectangle intersection graphs achieve $O(tn \cdot \log n / \log\log n)$.

\paragraph{Induced Tur\'an numbers and hereditary properties.}
Loh, Tait, and Timmons~\cite{Loh2016} showed that for any fixed $H$ one has $\text{ex}(n, \{H, K_{s,t}\text{-ind}\}) = O(n^{2-1/s})$.
Illingworth~\cite{Illingworth2021} asymptotically determined $\text{ex}(n, \{H, F\text{-ind}\})$ when $H$ is non-bipartite and $F$ is neither an independent set nor a complete bipartite graph, complementing the focus on complete bipartite forbidden induced subgraphs.
Nikiforov, Tait, and Timmons~\cite{Nikiforov2021} proved spectral strengthenings: if an induced $H$-free graph on $n$ vertices has spectral radius $\lambda(G) \geq Kn^{1-1/s}$ for appropriate $K$, then it contains $K_{s,t}$ as an induced subgraph.

The study of biclique-free graphs with hereditary constraints has seen significant recent progress.
Hunter, Milojevi\'c, Sudakov, and Tomon~\cite{Hunter2024} proved that for bipartite $H$ with maximum degree at most $k$ on one side, graphs avoiding induced $H$ and containing no $K_{s,s}$ have $O(s^{c} n^{2-1/k})$ edges, and conjectured that the dependence should match the ordinary Tur\'an number up to a constant factor; Axenovich and Zimmermann~\cite{AxenovichZimmermann2024} strengthened this result in bipartite host graphs, proving that for any $d \geq 2$ and any $K_{d,d}$-free bipartite graph $H$ where each vertex in one part has degree either at most $d$ or full degree, with at most $d-2$ full-degree vertices in that part, one has $\text{ex}(K_{n,n}, \{K_{t,t}, H\text{-ind}\}) = o(n^{2-1/d})$.
This result was built on insights of Janzer and Pohoata~\cite{JanzerPohoata2020}, who showed that $K_{k,k}$-free bipartite graphs with VC-dimension at most $d$ have $o(n^{2-1/d})$ edges, improving the KST bound when the VC-dimension is smaller than the biclique parameter.

\paragraph{VC-dimension, cographs, and the Erd\H{o}s-Hajnal conjecture.}
The \emph{VC-dimension} of a graph $G$ is defined as the VC-dimension of the hypergraph formed by closed neighborhoods $\{\mathscr{N}(v) : v \in V(G)\}$.
Intuitively, a graph has bounded VC-dimension when its neighborhood structure is not too complex---specifically, when no large vertex set can be ``shattered'' by neighborhoods.
Bousquet, Lagoutte, Li, Parreau, and Thomass\'e~\cite{Bousquet2014_identifyingCodes} proved a fundamental dichotomy: a hereditary graph class has finite VC-dimension if and only if it excludes some bipartite graph, some co-bipartite graph, and some split graph.
Remarkably, $P_4$ is simultaneously bipartite, co-bipartite, and split, so \emph{forbidding $P_4$ alone as an induced subgraph suffices to ensure bounded VC-dimension}. Indeed, all graphs $G$ such that the induced $G$-free graphs have bounded VC-dimension are induced subgraphs $G \subseteq P_4$.

Nguyen, Scott, and Seymour~\cite{Nguyen2023} resolved the Erd\H{o}s-Hajnal conjecture for graphs of bounded VC-dimension, proving that every $n$-vertex graph with VC-dimension at most $d$ contains a clique or stable set of size at least $n^{c_d}$ for some $c_d > 0$.
It was a well-known simple fact that cographs satisfy the Erd\H{o}s-Hajnal property with $c = 1/2$, meaning every $n$-vertex cograph contains a clique or independent set of size $\Omega(\sqrt{n})$~\cite{ErdosHajnal1989}.

The connection between induced subgraph avoidance and density bounds is highlighted by several fundamental results.
Fox and Sudakov~\cite{FoxSudakov2007} proved subexponential improvements to R\"odl's Theorem for hereditary graph classes, avoiding the tower-type bounds from the regularity Lemma and conjecturing that polynomial bounds should hold.
Gishboliner and Shapira~\cite{Gishboliner2023} showed that induced $P_4$-free graphs (cographs) have polynomial R\"odl bounds, exploiting the fact that $P_4$ is simultaneously bipartite, co-bipartite, and split.
Alon, Fischer, and Newman~\cite{Alon2007} proved that bipartite graph properties characterized by a finite set of forbidden induced subgraphs are efficiently testable with polynomial query complexity---specifically, one can $\varepsilon$-test membership in such a class using $\text{poly}(1/\varepsilon)$ queries, a dramatic improvement over the tower-type bounds that arise from the Szemer\'edi regularity Lemma.

\section{Biclique-Profiles} \label{sec:biclique_profiles}

For a graph $G$ let us define its \emph{biclique-sequence} by $\mathscr{S}(G)$, where for $s \in \NN$
\begin{align}
    \mathscr{S}_s(G) = \max \left\{ t \in \NN \ | \ K_{s,t} \subseteq G \right\}.
\end{align}
Here, we define $K_{0,t}$ as the empty graph on $t$ vertices, i.e., $\mathscr{S}_0(G) = |G|$. Observe that $\mathscr{S}_1(G)$ is the maximum degree.

Moreover, we define the maximum of the empty set as $-\infty$.
That is, the biclique-sequence of a graph is monotonically decreasing and constant $-\infty$ on indices greater than its vertex count. Furthermore, biclique-sequences are \emph{subdiagonal}, a property defined by the following two equivalent implications
\begin{alignat}{3}
    \mathscr{S}_{j}(G) < l
    &\;\implies\;& K_{j, l} \nsubseteq G
    &\;\implies\;& \mathscr{S}_{l}(G) < j \label{eq:subdiagonal_1},\\
    \mathscr{S}_{j}(G) \geq l
    &\;\implies\;& K_{j, l} \subseteq G
    &\;\implies\;& \mathscr{S}_{l}(G) \geq j \label{eq:subdiagonal_2}.
\end{alignat}

We might interpret $\mathscr{S}_{\infty}$ as the limit of the sequence, and deduce the following equivalence.
\begin{align}
    \mathscr{S}_{0}(G) < \infty &\iff \mathscr{S}_{\infty}(G) = -\infty \label{eq:finite_vertices_sequence}
\end{align}

\begin{Definition}[Biclique-profile] \label{def:biclique_profile}
A \emph{biclique-profile} is any decreasing subdiagonal sequence $\mathscr{P}$ in $\NN \cup \{ \infty, -\infty \}$.
\end{Definition}

Let us say a graph $G$ \emph{fulfills} a biclique-profile $\mathscr{P}$ in case that pointwise $\mathscr{S}(G) \leq \mathscr{P}$. Furthermore, a graph is $\mathscr{P}$-\emph{extremal} in case it fulfills $\mathscr{P}$ and it has the maximal edge count among all graphs with the same vertex count that fulfill $\mathscr{P}$. Let us say, a biclique $K_{s,t}$ is $\mathscr{P}$-\emph{problematic} if $\mathscr{P}_s < t$, i.e. every graph containing $K_{s,t}$ cannot fulfill $\mathscr{P}$.
We use the same notion for the extremal cograph problem.
For a biclique-profile $\mathscr{P}$ let us define its \emph{start-index} by the minimal $i \in \NN_{0}$ such that $\mathscr{P}_i < \infty$.

\subsection{Dynamic programming} \label{subsec:dynamic_programming}

The recursive structure of cographs via sum and product operations naturally leads to a dynamic programming approach for computing extremal cographs.
We first observe that biclique-sequences combine predictably under these operations.

\begin{Lemma}[Profile operations] \label{thm:profile_operations}
    Let $G_1, G_2$ be cographs on $n_1, n_2$ vertices respectively. Then:
    \begin{enumerate}[label=(\roman*)]
        \item For the sum $G_1 + G_2$:
        \begin{align}
            \mathscr{S}_s(G_1 + G_2) = \max\left(\mathscr{S}_s(G_1), \mathscr{S}_s(G_2)\right) \label{eq:profile_sum}
        \end{align}
        \item For the product $G_1 \times G_2$:
        \begin{align}
            \mathscr{S}_s(G_1 \times G_2) = \max_{a+c=s}\left(\mathscr{S}_a(G_1) + \mathscr{S}_c(G_2)\right) \label{eq:profile_product}
        \end{align}
    \end{enumerate}
\end{Lemma}

\begin{proofref}{Lemma}{thm:profile_operations}
    (i) In a sum $G_1 + G_2$, there are no edges between $V(G_1)$ and $V(G_2)$. Thus any $K_{s,t} \subseteq G_1 + G_2$ must lie entirely in $G_1$ or entirely in $G_2$, giving the pointwise maximum.

    (ii) In a product $G_1 \times G_2$, every vertex in $G_1$ is adjacent to every vertex in $G_2$. A $K_{s,t}$ subgraph can use $a$ vertices from one side of the biclique in $G_1$ and $c = s - a$ vertices from the other side, with the remaining vertices coming from $G_2$. The maximum over all such decompositions yields the max-convolution formula.
\end{proofref}

The key observation is that the number of edges also combines predictably:
\begin{align}
    \|G_1 + G_2\| &= \|G_1\| + \|G_2\| \label{eq:edges_sum} \\
    \|G_1 \times G_2\| &= \|G_1\| + \|G_2\| + n_1 \cdot n_2 \label{eq:edges_product}
\end{align}

For correctness of the algorithm presented below and for the proof of the Pumping Theorem we need the following technical helping Lemma.

\begin{Lemma}[Profile restriction] \label{thm:profile_restriction}
    Let $\mathscr{P}$ and $\mathscr{P}^{(1)}$ be biclique-profiles. Define $\mathscr{P}^{(2)}$ by
    \begin{align}
        \mathscr{P}^{(2)}_{c} := \min_{a \geq 0} \left( \mathscr{P}_{a+c} - \mathscr{P}^{(1)}_{a} \right) \label{eq:profile_restriction}
    \end{align}
    for all $c \in \NN_0$. Then $\mathscr{P}^{(2)}$ is a biclique-profile. Moreover, if $G = G_1 \times G_2$ and $G_1$ has biclique-sequence $\mathscr{S}(G_1) = \mathscr{P}^{(1)}$, then $G_2$ fulfills $\mathscr{P}^{(2)}$ if and only if $G$ fulfills $\mathscr{P}$.
\end{Lemma}

\begin{proofref}{Lemma}{thm:profile_restriction}
    We first verify that $\mathscr{P}^{(2)}$ is a valid biclique-profile by checking that it is monotonically decreasing and subdiagonal.

    \textbf{Monotonicity:} Let $c < c'$. Since $\mathscr{P}$ is monotonically decreasing, $\mathscr{P}_{a+c'} \leq \mathscr{P}_{a+c}$ for all $a \geq 0$. Thus,
    \begin{align*}
        \mathscr{P}^{(2)}_{c'} = \min_{a \geq 0} \left( \mathscr{P}_{a+c'} - \mathscr{P}^{(1)}_{a} \right) \leq \min_{a \geq 0} \left( \mathscr{P}_{a+c} - \mathscr{P}^{(1)}_{a} \right) = \mathscr{P}^{(2)}_{c}.
    \end{align*}

    \textbf{Subdiagonality:} Suppose $\mathscr{P}^{(2)}_{j} < l$ for some $j, l \in \NN_0$. We need to show $\mathscr{P}^{(2)}_{l} < j$. By definition, there exists $a^* \geq 0$ such that $\mathscr{P}_{a^*+j} - \mathscr{P}^{(1)}_{a^*} < l$. This implies $\mathscr{P}_{a^*+j} < l + \mathscr{P}^{(1)}_{a^*}$.

    Let us set $a' \coloneq \mathscr{P}^{(1)}_{a^*}$. By the subdiagonality of $\mathscr{P}$ applied with indices $a^* + j$ and $l + a'$, we obtain $\mathscr{P}_{l + a'} < a^* + j$.

    Now, by subdiagonality of $\mathscr{P}^{(1)}$, since $a' \geq a'$, we have $\mathscr{P}^{(1)}_{a'} \geq a^*$. In the minimum~\eqref{eq:profile_restriction}:
    \begin{align*}
        \mathscr{P}^{(2)}_{l} \leq \mathscr{P}_{l + a'} - \mathscr{P}^{(1)}_{a'} < (a^* + j) - a^* = j.
    \end{align*}

    \textbf{Fulfillment:} Now assume $G = G_1 \times G_2$ fulfills $\mathscr{P}$ and $\mathscr{S}(G_1) \leq \mathscr{P}_1$. By Lemma~\ref{thm:profile_operations}~(ii):
    \begin{align*}
        \mathscr{S}_s(G) = \max_{a+c=s} \left( \mathscr{S}_a(G_1) + \mathscr{S}_c(G_2) \right).
    \end{align*}
    Since $G$ fulfills $\mathscr{P}$, for all $s \geq 0$: $\mathscr{S}_s(G) \leq \mathscr{P}_s$. Thus for any $a, c$ with $a + c = s$:
    \begin{align*}
        \mathscr{P}^{(1)}_a + \mathscr{S}_c(G_2) \leq \mathscr{P}_s = \mathscr{P}_{a+c}.
    \end{align*}
    Taking the minimum over all $a \geq 0$ (letting $s = a + c$ range appropriately):
    \begin{align*}
        \mathscr{S}_c(G_2) \leq \min_{a \geq 0} \left( \mathscr{P}_{a+c} - \mathscr{P}^{(1)}_{a} \right) = \mathscr{P}^{(2)}_{c}.
    \end{align*}
    Hence $G_2$ fulfills $\mathscr{P}^{(2)}$. For the other direction let us assume that $G_2$ fullfills $\mathscr{P}^{(2)}$.
    By Lemma~\ref{thm:profile_operations}~(ii):
    \begin{align}
        \mathscr{S}_c
            \leq \max_{a+c = s} \left( \mathscr{P}^{(1)}_a + \mathscr{P}^{(2)}_c \right)
            \leq \max_{a+c = s} \left( \mathscr{P}^{(1)}_a + \mathscr{P}_{a + c} - \mathscr{P}^{(1)}_a \right)
            = \mathscr{P}_{s}.
    \end{align}
    This closes the proof of Lemma \ref{thm:profile_restriction}.
\end{proofref}

\begin{Definition}[Pareto frontier] \label{def:pareto_frontier}
    For a set of pairs $\{(\mathscr{P}_j, e_j)\}_{j \in J}$ consisting of biclique-profiles and edge counts, the \emph{Pareto frontier} is the set of pairs $(\mathscr{P}_j, e_j)$ such that there is no other pair $(\mathscr{P}_k, e_k)$ with $\mathscr{P}_k \leq \mathscr{P}_j$ (coordinatewise) and $e_k \geq e_j$, with at least one strict inequality.
\end{Definition}

The algorithm proceeds by dynamic programming over the vertex count $n$. For each $n$, we maintain a registry $\mathcal{R}_n$ of pairs $(\mathscr{P}, \mathcal{G}_{\mathscr{P}})$, where $\mathcal{G}_{\mathscr{P}}$ is the set of cographs on $n$ vertices having profile $\mathscr{P}$ with maximum edge count among cographs with that profile.

\begin{center}
\fbox{\parbox{0.9\textwidth}{
\textbf{Algorithm: Dynamic Programming for Extremal Cographs} \label{alg:dynamic_programming}
\begin{enumerate}[label=\arabic*.]
    \item \textbf{Base case.} Set $\mathcal{R}_1 \leftarrow \{((1, 0, \ldots), \{K_1\})\}$.

    \item \textbf{Inductive step.} For each $n \geq 2$:
    \begin{enumerate}[label=(\alph*)]
        \item Initialize candidate set $\mathcal{C}_n \leftarrow \emptyset$.
        \item For each partition $n = n_1 + n_2$ with $1 \leq n_1 \leq n_2$:
        \begin{enumerate}[label=(\roman*)]
            \item For each $(\mathscr{P}_1, \mathcal{G}_1) \in \mathcal{R}_{n_1}$ and $(\mathscr{P}_2, \mathcal{G}_2) \in \mathcal{R}_{n_2}$:
            \begin{itemize}
                \item Compute $\mathscr{P}_{+} = \mathscr{S}(\mathscr{P}_1 + \mathscr{P}_2)$ via~\eqref{eq:profile_sum}.
                \item Compute $\mathscr{P}_{\times} = \mathscr{S}(\mathscr{P}_1 \times \mathscr{P}_2)$ via~\eqref{eq:profile_product}.
                \item Compute edge counts $e_{+}$, $e_{\times}$ via~\eqref{eq:edges_sum}--\eqref{eq:edges_product}.
                \item Add $(\mathscr{P}_{+}, e_{+}, G_1 + G_2)$ and $(\mathscr{P}_{\times}, e_{\times}, G_1 \times G_2)$ to $\mathcal{C}_n$ for each $G_1 \in \mathcal{G}_1$, $G_2 \in \mathcal{G}_2$.
            \end{itemize}
        \end{enumerate}
        \item \textbf{Pareto filtering.} Reduce $\mathcal{C}_n$ to the Pareto frontier with respect to the partial order $(\mathscr{P}, e) \preceq (\mathscr{P}', e')$ iff $\mathscr{P} \geq \mathscr{P}'$ and $e \leq e'$.
        \item Group by profile: $\mathcal{R}_n \leftarrow \{(\mathscr{P}, \{G : (\mathscr{P}, e, G) \in \mathcal{C}_n, e = e_{\max}(\mathscr{P})\})\}$.
    \end{enumerate}

    \item \textbf{Query.} Given a biclique-profile $\mathscr{P}$ with $\mathscr{P}_0 = n$, return cographs from $\mathcal{R}_n$ whose profiles are $\leq \mathscr{P}$ (coordinatewise) with maximum edge count.
\end{enumerate}
}}
\end{center}

The lattice-based optimization (step 2c) prunes dominated profile pairs at the profile level \emph{before} expanding to individual cograph pairs, reducing the combinatorial explosion.

For practical computation with a pruning constraint $K_{s,t}$, we can truncate profiles to length $s+1$ since only the first $s$ entries determine whether a cograph fulfills the constraint $\mathscr{P}_s < t$. Moreover, when only interested in $(s,t)$-extremal cographs we can prune all profiles $\mathscr{P}$ where $\mathscr{P}_s \geq t$.

It is easy to see that our approach is correct: By Lemma \ref{thm:profile_restriction} in any $(s,t)$-extremal cograph all components are extremal with respect to a biclique profile that depends on the remaining cograph. Thus, by combining all cographs extremal to suitable profiles, the search considers all relevant candidates for extremal cographs.

By step (c) we lose the ability to enumerate all extremal cographs. We remark that in order to enumerate all extremal cographs one could simply skip step (c) - at the cost of significantly more memory consumption and run time. However, it is easy to verify that also with step (c) one still finds at least one extremal cograph for each $n \in \NN$.

\subsection{The Lifting Theorem} \label{subsec:lifting}

We say a biclique-profile $\mathscr{P}$ \emph{dominates} a distinct profile $\mathscr{P}'$ if $\mathscr{P}_i \geq \mathscr{P}'_i$ for all $i \in \NN_0$ and there is at least one inequality, and write $\mathscr{P} > \mathscr{P}'$.
This defines a partial order on profiles, which forms a lattice structure.

In the following Theorem we make use of the so-called ``lifting'' procedure, where we replace the sum of components in an extremal structure by a product of a small core and the sum of smaller structures.

\begin{Theorem}[Lifting Theorem] \label{thm:lifting}
    Let $\mathscr{P}$ be a biclique-profile with start-index $s$ such that $\mathscr{P}_s = t-1$. For any $\mathscr{P}$-extremal cograph $G$ there is a finite index set $I$ and a mapping $\sigma: I \to [t-1]$ such that:
    \begin{enumerate}[label=(\roman*)]
        \item $G$ admits the decomposition:
        \begin{align}
            G = \sum_{i \in I} G_{i,1} \times G_{i,2}
        \end{align}
        where $|G_{i,1}| = \sigma(i)$ and $|G_{i,2}| \geq \sigma(i)$ for all $i \in I$.
        \item For any $1 \leq j \leq s-1$, the fiber $|\sigma^{-1}(j)| \leq 1$.
    \end{enumerate}
\end{Theorem}

\begin{Observation}[Component size bound] \label{obs:component_size_bound}
    Let $\mathscr{P}$ be a biclique-profile with start-index $s$ and let $G$ be a connected cograph fulfilling $\mathscr{P}$ such that $G = G_1 \times G_2$ with $|G_1|, |G_2| \geq s$. Then $|G| \leq 2\mathscr{P}_s$.
\end{Observation}

\begin{proofref}{Observation}{obs:component_size_bound}
    It is clear that both $|G_1| \leq \mathscr{P}_s$ and $|G_2| \leq \mathscr{P}_s$.
\end{proofref}

\begin{proofref}{Theorem}{thm:lifting}   
    Let $G$ be a $\mathscr{P}$-extremal cograph.
    Each of its connectivity component $C_i$ for $i \in I$ is a product. We write $C_i = G_{i,1} \times G_{i,2}$ and choose the decomposition such that $\sigma(i) \coloneqq |G_{i,1}|$ is minimal. Note that $\sigma(i) < t$ must hold. If not, $|G_{i,1}| \geq t$ and $|G_{i,2}| \geq t$, implying $K_{t,t} \subseteq C_i$. By (\ref{eq:subdiagonal_1}) and the start-index property of $s$ we have $t-1 \geq s-1$, so $s \leq t$. This implies $K_{s,t} \subseteq G$, which violates $\mathscr{P}_s = t-1$.
    
    To prove (ii), assume for a contradiction that there exist distinct indices $i_1, i_2 \in I$ such that $\sigma(i_1) = \sigma(i_2) = j$ for some $j \in [s-1]$. Now we continue by a lifting argument.
    
    Consider the cograph $G'$ formed by replacing the components $C_{i_1} + C_{i_2}$ with a single component:
    \begin{align*}
        C' \coloneqq E_j \times (G_{i_1,2} + G_{i_2,2} + K_j).
    \end{align*}
    This cograph $G'$ has the same vertex count as $G$. We claim $G'$ still fulfills $\mathscr{P}$.

    Let $s' \geq s$ and $p \coloneqq \mathscr{P}_{s'}$.
    It is clear that $p \geq s-1$ as otherwise by (\ref{eq:subdiagonal_1}) $\mathscr{P}_{s-1} < \infty$, a contradiction to the definition of $s$ as the start-index. Hence any biclique that would violate $\mathscr{P}$ has size at least $2s$.
    
    It is clear that $G'$ does not contain $K_{s', p+1}$ unless the product $C'$ does contain it. However, this could only be the case if $E_j \times K_j$ does contain it, which is not possible as $2s > 2 j$.
    
    Furthermore, $G'$ has strictly more edges than the sum $C_{i_1} + C_{i_2}$.
    Indeed, the cross edges between $G_{i_1, 1}$ and $G_{i_2, 2}$ are covered by the edges between $E_j$ and $G_{i_1, 2}$ as well as $E_j$ and $G_{i_2, 2}$.
    Moreover, $\| G_{i_1, 1} \| \leq K_j$ and the remaining cross edges exceed the edges inside $\| G_{i_2, 1} \|$ as $\| G_{i_2, 1} \| \leq \binom{j}{2} < j^2 = \| E_j, K_j \|$.
    This contradicts the edge-maximality of $G$.
\end{proofref}

\section{Proofs of the Main Results} \label{sec:main_results}

\subsection{Special cases}

Before we consider the general case, we prove the special cases stated in Section~\ref{sec:main_results_intro}.
It is clear that $(t-1)$-extremal cographs are cographs of maximal degree $t-1$.
Before we study the extremal cographs, let us study regular cographs.

\begin{Lemma} \label{thm:regular_cographs}
    Let $n \in \mathbb{N}$.
    \begin{enumerate}
        \item If $n$ is even, there exists a $d$-regular cograph on $n$ vertices for every $0 \le d < n$.
        \item If $n$ is odd, for any $0 \leq d < n$ such that $d$ is even and $2d \neq n-1$ there exists a $d$-regular cograph on $n$ vertices.
    \end{enumerate}
\end{Lemma}

\begin{proofref}{Lemma}{thm:regular_cographs}
    
We proceed by strong induction on $n$.
The base cases $n \in \{1, 2, 3 \}$ are trivial. Let us assume that $n \geq 4$ and the Theorem holds for all integers $k < n$.

\paragraph{case 1: $n$ is even.}
We need to construct a cograph for any $0 \le d < n$.

\paragraph{subcase 1.1: $d$ is odd.}
Consider the target degree $d' = n - 1 - d$. Since $n$ is even and $d$ is odd, $d'$ is even. Also $0 \le d' < n$. By the arguments in Subcase 1.2 (below), there exists a $d'$-regular cograph $G'$. The complement $\overline{G'}$ is a cograph with degree $n - 1 - d' = d$.

\paragraph{subcase 1.2: $d$ is even.}
We work with two building blocks: $K_{d+1}$ and the $d$-regular complete multipartite  cograph
\begin{align}
    H_d \coloneq \prod_{j \in \left[ \frac{d}{2} + 1 \right]} E_2.
\end{align}

\paragraph{subsubcase 1.2.1: $3d \neq n-2$.}
    We construct $G = K_{d+1} + G'$.
    We need $G'$ to be a $d$-regular cograph on the remaining $n' = n - (d+1)$ vertices.
    Since $n$ is even and $d+1$ is odd, it is clear that $n'$ is odd.
    Moreover, 
    \begin{align}
        3d \neq n - 2 \implies 2d \neq (n - d - 1) - 1 = 2n' - 1.
    \end{align}
    Thus, induction yields $G'$.

\paragraph{subsubcase 1.2.2: $3d = n-2$.}
    As the previous construction fails, we construct $G = H_d + G$. 
    We need $G$ to be a $d$-regular cograph on $n' = n - (d+2)$ vertices.
    Since $n$ and $d$ are even, $n'$ is even.
    Either $d = 0$ and $n = 2$ so $G = H_d$ or it is obvious that $2d < n-2$ which implies that $d < n - d - 2 = n'$. In the latter case induction yields the required $G'$.

\paragraph{case 2: $n$ is odd}
We need to construct a cograph for even $d$ with $0 \leq d < n$ satisfying $2d \neq n-1$.

\paragraph{subcase 2.1: $d < (n-1)/2$.}
    We construct $G = K_{d+1} + G'$.
    $G'$ must be on $n' = n - (d+1)$ vertices. Since $n$ is odd and $d+1$ is odd, $n'$ is even. As $2d < n-1$ we have $d < n-d-1=n'$ and induction yields $G'$.

\paragraph{subcase 2.2: $d > (n-1)/2$.}
    Consider the complement degree $d' = n - 1 - d$.
    Apparently $2d' = 2n - 2 - 2d < 2n - 2 - (n-1) = n - 1$, so subcase 2.1 yields $d'$-regular cograph $G'$ on $n$ vertices. It's complement is $d$-regular.

    This closes the proof of Theorem \ref{thm:regular_cographs}.
\end{proofref}

Let $t,n \in \NN$ with $n < t$. It is clear that any $(1,t)$-extremal cograph on $n$ vertices is a clique.

\beginProofOfTheoremClickable{thm:pumping_stars}
    In case $n$ is even or $t$ is even and $2(t-1) \neq n-1$ the claim is evident from Lemma \ref{thm:regular_cographs}, so we may assume that $n$ is odd. Let $G$ be $(1,t)$-extremal on $n$ vertices.
    Let us assume for a contradiction that there are several connectivity components $C_1, \dots, C_l$ that are not $(t-1)$-regular and a $(t-1)$-regular cograph $G'$ such that $G = \sum_{1 \leq j < l} C_j + G'$. We may assume that $l \geq 2$ is minimal among all $(1,t)$-extremal cographs.

    \paragraph{case 1: $\sum_{1 \leq j \leq l} |C_j|$ is even.}
        In case $\sum_{1 \leq j \leq l} |C_j| \leq t-1$, one could replace $\sum_{1 \leq j \leq l} C_j$ by a clique, otherwise by a $(t-1)$-regular cograph. A contradiction to edge-maximality in both cases.

    \paragraph{case $\sum_{1 \leq j \leq l} |C_j|$ is odd.}
        Then $|G'|$ is even. Either $\sum_{1 \leq j < l} |C_j|$ is even or $|C_l|$ is even, we may assume the latter. Then however, also $|C_l + G'|$ is even and if $|C_l + G'| \leq t-1$ we could replace it by a clique, otherwise by a $(t-1)$-regular cograph. Both are contradictions to our structural assumption.

    It is clear that any connected $(1,t)$-extremal cograph has size at most $2(t-1)$. In case it has size $2(t-1)$ it however is $K_{t-1, t-1}$ - a $(t-1)$-regular cograph. This shows that $|C_1| \leq 2t-3$.
\end{proof}

We want to remark that for example all complete multipartite graphs with equally sized partition classes are connected regular cographs. However, there are many more, as there are even connected regular cographs with different vertex degrees in different multiplication components---consider the $4$-regular cograph $2K_2 \times E_3$.

\beginProofOfTheoremClickable{thm:classification_2_t}
    Let $t \in \{ 2,3 \}$ and $n \in \NN$ and let $G$ be a $(2,t)$-extremal cograph.
    
    Let us show that any component $H$ has a complete vertex $v$. Then by the Lifting Theorem there can only be one connectivity component and $G$ is the product of $v$ with $G'$, where it is apparent that $G'$ is $(1,t)$-extremal.
    
    In case that $H$ is not a complete vertex it is the product of two cographs $H_{1}$ and $H_{2}$ on $n_1$ and $n_2$ vertices respectively.

    \paragraph{case $t = 2$.}
        Let us assume for a contradiction that $H_j$ does not contain a complete vertex. But then $n_1 \geq 2$ and $n_2 \geq 2$, a contradiction as $K_{2,2} \subseteq H$.

    \paragraph{case $t = 3$.}
        We may assume that $n_1 \geq 2$ and $n_2 \geq 2$.

        If $n_1 \geq t$ it is clear that $n_2 = 1$ and we have found our extremal vertex. By symmetry we may assume that $n_1 \leq 2$ and $n_2 \leq 2$.
        Thus $|H| = 4 < |K_{2,3}|$ so $H$ is complete and we can pick an arbitrary complete vertex.
\end{proof}

Let us prove a technical Lemma that serves as a induction base in the Pumping Theorem for biclique-profiles \ref{thm:pumping_profile}.

\begin{Lemma} \label{thm:classification_2_t_function}
    For any $t \geq 2$
    \begin{align}
        \text{ex} \left( n, \{ K_{2,t}, P_4 - \text{ind} \} \right) < \frac{t+1}{2} \cdot n
    \end{align}
\end{Lemma}

\beginProofOfLemmaClickable{thm:classification_2_t_function}
    Let us prove the Theorem by induction on $t$. For the base $t=2$ we know that any $(2,2)$-extremal cographs have a complete vertex. All the other vertices have degree at most $t-1$.
    Hence $\text{ex} \left( n, \{ K_{2,t}, P_4 - \text{ind} \} \right) \leq \left(1 + \frac{t-1}{2} \right) \cdot (n-1) < \frac{t+1}{2} \cdot n$.

    For the step, let us assume that $t \geq 3$ and $n \in \NN$.
    Let $H_1, \dots, H_l$ be the connectivity components of a $(2,t)$-extremal cograph.
    In case, $H_1$ contains a complete vertex by a similar argument as above it follows that $\| H_1 \| < \frac{t+1}{2} \cdot |H_1|$.
    Otherwise, we know that there are cographs $H_{1,1}$ and $H_{1,2}$ such that $H_1 = H_{1,1} \times H_{1,2}$ on $n_1$ and $n_2$ vertices respectively and may assume that $2 \leq n_1 \leq n_2 \leq t-1$.

    In case that $n_1 = t-1$ we know that there is no $K_{1,2}$ in $H_{1,2}$ and $H_{1,2}$ is a matching. This is, there are at most $\frac{n_2}{1}$ edges in $H_{1,2}$.
    In case that $n_1 < t-1$, we know that there is no $K_{2, t-n_1}$ in $H_{1,2}$ so induction yields that $\| H_{1,2} \| < \frac{t-n_1 + 1}{2} \cdot n_2$.
    The same holds for $n_2$ and $H_{1,1}$ by symmetry.
    Hence, in case $n_1=n_2=t-1$
    \begin{align}
        \| H \|
            \leq (t-1)^2 + \frac{t-1}{2} + \frac{t-1}{2}
            < \frac{t+1}{2} \cdot 2(t-1)
            = \frac{t+1}{2} \cdot |H_1|.
    \end{align}
    In case $n_1<n_2=t-1$
    \begin{align}
        \| H_1 \|
            \leq n_1 \cdot (t-1) + \frac{n_1}{2} + \frac{t - n_1 + 1}{2} \cdot (t-1)
            \leq n_1 \cdot \frac{t}{2} + \frac{t+1}{2} \cdot n_2
            < \frac{t+1}{2} \cdot |H_1|.
    \end{align}
    Lastly, in case $n_1, n_2 < t-1$ we have
    \begin{align}
        \| H_1 \|
            < n_1 \cdot n_2 + \frac{t - n_2 + 1}{2} \cdot n_1 + \frac{t - n_1 + 1}{2} \cdot n_2.
            = \frac{t+1}{2} \cdot |H_1|
    \end{align}

    Combining similar inequalities for the other connectivity components, this proves the claim and closes the proof of Lemma \ref{thm:classification_2_t_function}.
\end{proof}

\beginProofOfTheoremClickable{thm:classification_3_3}
    Let us prove the claim by induction on $n$.
    In case that $2 \leq n < 6$ it is clear that $G$ is a clique and therefore the product of an edge with another cograph. This solves the base case.
    
    \textbf{step.} We might assume that $n \geq 6$.
        First, let us prove that $G$ is connected. Assume for a contradiction that $G = H_1 + H_2$ for two cographs $H_1$ and $H_2$.
        Then, either $H_1 \geq 2$ and induction yields that we might replace $H_1$ by an edge-maximal cograph not containing $K_{3,3}$ that has a complete vertex $v_1$. Otherwise, $H_1$ itself is a singleton and a complete vertex itself. The same argument yields a complete vertex $v_2$ of $H_2$. However,
        \begin{align}
            \tilde{G} \coloneqq v_1 \times \left( \left( H_1 - v_1 \right) + \left( H_2 - v_2 \right) + v_2 \right)
        \end{align}
        has one more edge than $G$ and also does not contain $K_{3,3}$, a contradiction.

        Hence, we know that $G = H_1 \times H_2$ for two cographs $H_1$ and $H_2$ on $n_1$ and $n_2$ vertices respectively.
        Let us assume that $n_1 \leq n_2$. Moreover, we might assume that our choices of $H_1$ and $H_2$ maximize $n_1$ under this restriction.
        It is clear that $n_1 < 3$, so to show $n_1 = 2$ we only need to show that $n_1 \neq 1$.
        
        Let us assume for a contradiction that $H_1$ is made of a single vertex $v$. Since we assumed $H_1$ to be maximal under the restriction $n_1 \leq n_2$, we might assume that $H_2$ is disconnected. Indeed, otherwise it is a product of two components, and in case their size differs, we could consider the smaller of the two together with $v$ as $H_1$ or in case their size are equal, we could consider one of them as $H_2$.
        This is, there is $l \geq 2$ and cographs $H_{2,1},\dots,H_{2,l}$ such that $H_2 = \sum_{j \in [l]} H_{2,j}$. However, it is easy to see that each $H_{2,j}$ is $(2,3)$-extremal so Theorem \ref{thm:classification_2_t}
        yields a complete vertex $v_{j}$ in $H_{2,j}$.
        However,
        \begin{align}
            \tilde{G} \coloneqq v \times v_1 \times \left( \sum_{j \in [l]} \left( H_j - v_j \right) + \sum_{j \in [l]\setminus\{1\}} v_j \right)
        \end{align}
        does not contain $K_{3,3}$ and has $l-1$ more edges than $G$, a contradiction.

        Thus, we know that $H_1$ does contain two vertices. It is left to show that there is an extremal cograph $G$ with two complete vertices, i.e. the respective $H_1$ is an edge.
        It is clear that $H_2$ does not contain a $K_{1,3}$.
        
        \begin{Claim} \label{claim:13}
            Any edge-maximal cograph $G$ not containing $K_{1,3}$ is the sum of cliques of size at most $3$ and four-cycles.
        \end{Claim}

        \begin{proofref}{Claim}{claim:13}
            It is easy to see that any connectivity component of $G$ can have at most $4$ vertices. It is easy to check manually that indeed, for any vertex count of at most $4$ any edge-maximal and connected cograph not containing $K_{1,3}$ is either a clique or a four-cycle.
            This concludes the proof of Claim \ref{claim:13}.
        \end{proofref}

        Claim \ref{claim:13} yields that $H_2$ is the sum of cliques of size at most $3$ and four-cycles.
        We might assume that it is the sum of $a$ triangles, $b$ four-cycles, and a clique of size $s \in \{ 0, 1, 2 \}$.
        Then,
        \begin{align}
            \| H_2 \| = 3 \cdot a + 4 \cdot b + \mathbb{1}\{ s = 2 \}
        \end{align}
        In case that $H_2$ does not contain a four-cycle, by edge-maximality it is clear that $H_1$ is an edge.
        In case that $H_2$ does contain a four-cycle, $H_1$ cannot be an edge, as otherwise a $K_{3,3}$ would arise.
        Let us check all cases to see that in case $H_2$ does contain a four-cycle, we can replace the four-cycles with cliques and $H_1$ by an edge and obtain a $K_{3,3}$-free cograph of the same edge count.

        \paragraph{case $n_2 \bmod 3 = 0$.} In this case, $H_2$ can be chosen to be made out of triangles and $H_1$ would be an edge.

        \paragraph{case $n_2 \bmod 3 = 1$.} In this case $H_2$ can be chosen to be made of triangles and either a singleton or a single four-cycle. In the first case, $H_1$ is an edge. In the latter case we can replace the four-cycle by a triangle and a singleton and add the edge in $H_1$ for compensation.

        \paragraph{case $n_2 \bmod 3 = 2$.}  In this case $H_2$ can be chosen to be made of triangles and either an edge or two four-cycles. In the first case, $H_1$ is an edge. In the latter case, we can replace the two four-cycles by two triangles and an edge and add the edge in $H_1$ for compensation.
        This completes the proof of Theorem \ref{thm:classification_3_3}.
\end{proof}

\subsection{The general Pumping Theorem}

To proof our Main Theorem we use the following well-known result.

\begin{Theorem}[Davenport] \label{thm:davenport}
    For any $n \in \NN$, let $A$ be a sequence of $n$ integers. There exists a non-empty subsequence of $A$ whose sum is divisible by $n$.
\end{Theorem}

We first introduce the cotree representation.
The \emph{cotree} $T_G$ of a cograph $G$ is a rooted labeled tree representing its construction sequence.
Inner vertices are labeled by $+$ or $\times$ for sum and product respectively (each of arity at least two), while leaves correspond to vertices of $G$.
The root operation corresponds to the final step in the construction.
Note that $G$ is connected if and only if the root is labeled $\times$.

\begin{figure}[ht]
    \centering
    \includegraphics[width=0.5\textwidth]{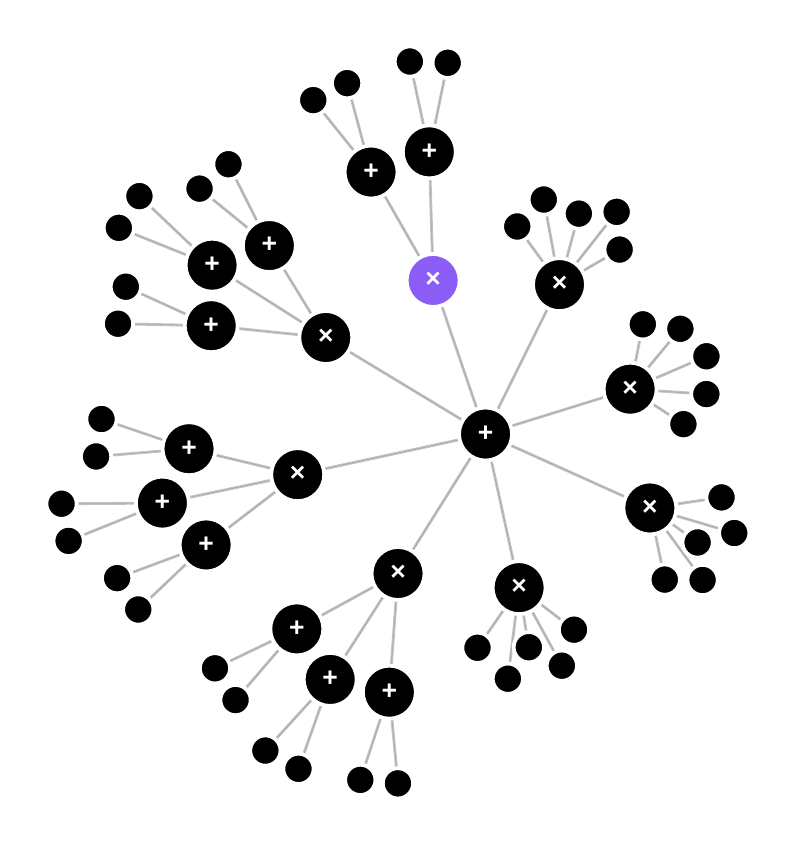}
    \caption{Cotree of a $K_{5,5}$-free extremal cograph on 42 vertices. The root is marked as purple. The corresponding cograph is the product of $K_{2,2}$ together with the sum of four $5$-cliques and three complete multipartite $K_{2,2,2}$.}
    \label{fig:cotree_example}
\end{figure}

For uniqueness, we require cotrees to be \emph{reduced}: no two adjacent inner vertices share the same label (using higher-arity operations instead).
The \emph{height} of a cotree $\text{height}(T_G)$ is the maximum number of edges on any root-to-leaf path.
For convenience, we write $\text{height}(G)$ directly.

The clique number behaves naturally with respect to cograph operations: $\omega(G_1 + G_2) = \max(\omega(G_1), \omega(G_2))$ and $\omega(G_1 \times G_2) = \omega(G_1) + \omega(G_2)$. Indeed, algebraically it can be seen as a valuation into the max-plus algebra.
This gives the bound
\begin{align}
    \text{height}(G) \leq 2 \cdot \omega(G) + 1. \label{eq:height_omega}
\end{align}
In particular, any $(s,t)$-extremal cograph has height at most $2(s+t) - 1$.

For $u \in V(T)$ let us denote the subtree of $u$ and its descendants with respect to the root of $T$ by $T(u)$. Let us transfer this notation to $G$ via
\begin{align}
    G(u) \coloneqq G[\mathscr{L}(u)],
\end{align}
where $\mathscr{L}(u)$ are the cograph vertices that correspond to the leaves of $T(u)$.
It is easy to see that $G(u)$ is the cograph corresponding to the cotree $T(u)$.

Before we prove the main result, we prove a Pumping Theorem for $\mathscr{P}$-extremal cographs for some biclique-profile $\mathscr{P}$. For this let us introduce the notation
\begin{align}
    \text{ex} \left( n, \{ \mathscr{P}, P_4 - \text{ind} \} \right)
\end{align}
to express the edge-count of $\mathscr{P}$-extremal cographs on $n$ vertices.
We remark that $\mathscr{P}$-extremal cographs are exactly the edge-maximal cographs not containing any of the $\mathscr{P}$-problematic bicliques.

\begin{Theorem}[Pumping Theorem for biclique-profiles] \label{thm:pumping_profile}
    Let $s, t \in \NN$ and $\mathscr{P}$ be a biclique-profile with start-index $s$ and start value $\mathscr{P}_s = t-1$ as well as non-negative limit.
    There exist $\alpha \in \mathbb{Q}_{\geq0}$, $R, N_{s,t} \in \NN$, and core graphs $G_0, \ldots, G_{R-1}$ each of size at most $N_{s,t}$ with designated components $H_j \subseteq G_j$ where any vertex in $H_j$ has exactly the same neighborhood outside of $H_j$ such that for $n$ large enough, there is a $\mathscr{P}$-extremal cograph on $n$ vertices obtained by pumping $H_{(n \bmod R)}$ inside of $G_{(n \bmod R)}$.
    There are constants constants $a_0, \ldots, a_{R-1} \in \mathbb{Q_{\leq 0}}$ such that for large $n$
    \begin{align}
        \text{ex}(n, \{\mathscr{P}, P_4\text{-ind}\}) = a_{(n \bmod R)} + \alpha n
    \end{align}
    Moreover, $\alpha \leq s - 1 + \frac{t-1}{2}$ and in case $s \geq 2$ we have the upper bound $\text{ex}(n, \{\mathscr{P}, P_4\text{-ind}\}) < \left( s-1 + \frac{t-1}{2} \right) \cdot n$ for any $n \in \NN$.
    In case that $H_q$ is a connectivity component in $G_q$ for some $q \in \mathbb{Z} / \mathbb{Z}$ we certainly know that $\alpha < s - 1 + \frac{t-1}{2}$.
\end{Theorem}

\paragraph{Proof outline.}
The proof of the Pumping Theorem for biclique-profiles proceeds in six steps:
\begin{enumerate}
    \item \textbf{Bounded degree of product vertices:} We show that product vertices in the cotree have bounded degree, since too many children would create a forbidden biclique.
    \item \textbf{Existence of high-degree addition vertices:} For large enough cographs, there must exist addition vertices with many children, as otherwise the cotree height bound would limit the vertex count.
    \item \textbf{Structure analysis:} We analyze the structure of subgraphs rooted at high-degree addition vertices, using the Lifting Theorem to constrain which children can be large.
    \item \textbf{Reduction to core with pumped component:} Using the pigeonhole principle and edge contribution arguments, we show that most small components can be replaced by copies of a single ``pumping component'' without losing edges, reducing to a bounded-size core graph.
    \item \textbf{Cyclic repetition:} Since there are finitely many core graphs and pumping components, the extremal function eventually becomes periodic in $n$, with period $R$ dividing the pumping component sizes.
    \item \textbf{Upper bounds:} We establish the claimed bounds on $\alpha$ and the constants $a_0, \ldots, a_{R-1}$ by induction on $(s, t-1)$ in lexicographic order.
\end{enumerate}
This result is then applied in Theorem~\ref{thm:pumping} to the specific biclique-profile forbidding $K_{s,t}$, where matching lower bounds from explicit constructions determine the exact linear coefficient $\alpha = s - 1 + \frac{t-1}{2}$.

\beginProofOfTheoremClickable{thm:pumping_profile}
    In case $s=1$, the proof of the claim is simpler but instructive to elaborate as it already demonstrates arguments used in the general case.
    By Observation \ref{obs:component_size_bound} any connectivity component of any $\mathscr{P}$-extremal cograph $G$ has size at most $2(t-1)$.
    In any $\mathscr{P}$-extremal graph $G$ there is a component $H$ with highest edge density.
    Using the Davenport Theorem \ref{thm:davenport} for any set of $|H|$ many connectivity components, one finds a subset such that the sum of their vertex count is divisible by $|H|$.
    Hence, one can replace this subset by pumping up $H$ appropriately often. Iterating this, one obtains an extremal graph $\hat{G}$ that contains at most $|H| - 1 \leq 2(t-1) - 1$ connectivity components not isomorphic to $H$ and many other components isomorphic to $H$.
    Deleting all $H$-isomorphic subgraphs from $\hat{H}$ one obtains a core graph $G_{\text{core}}$ on at most $2(t-1)^2$ vertices.
    The set of all such core graphs obtained by this procedure on all $\mathscr{P}$-extremal cographs is finite.
    Hence, there are candidates where pumping up the densest component $H$ yields extremal graphs for infinitely many $n$.
    Hence, one may choose a period $R$ such that for large $n$ one obtains $\mathscr{P}$-extremal graphs by pumping up one component in a core graph depending on $n \bmod R$.
    As we see in the general case, the density of this pumping components needs to be a constant $\alpha$, giving rise to the form of the extremal function $\text{ex}(n, \{\mathscr{P}, P_4\text{-ind}\}) = a_{(n \bmod R)} + \alpha n$.
    We see in the general case an argument why the constants $a_j$ are non-positive. As the maximum degree in any $\mathscr{P}$-extremal cograph is $t-1$ it is clear that $\alpha < \frac{t-1}{2}$.
    
    Now, let us assume that $s \geq 2$.
    Let us fix $D \coloneqq 4t^2 + s$ and $N \coloneqq D ^ {2(s+t) + 1} + 1$.
    Let $G$ be a $\mathscr{P}$-extremal cograph on $n \geq N$ vertices. Let $T_G$ be its reduced cotree.
    
    \paragraph{1. Bounded degree of product vertices:}
    Let $u$ be any product vertex in $T_G$. Let $v_1, \dots, v_k$ be its children. Since $G$ is a cograph, the subgraph induced by $u$ contains a complete multipartite graph with parts of sizes $|G(v_1)|, \dots, |G(v_k)|$. Since every child in the cotree has at least one vertex, if $k \geq s+t$, one could find a $K_{s,t}$ in $G$. Thus, any product vertex has strictly less than $s+t \leq D$ children.
    
    \paragraph{2. Existence of a high degree addition vertex:}
    Let us assume for a contradiction that all addition vertices have at most $D$ children. By (\ref{eq:height_omega}) it is clear that $\text{height}(G) \leq 2(s+t) + 1$. This implies however that $n \leq D ^ {2(s+t) + 1}$, a contradiction.
    Thus, we know that there exists a vertex $v$ with at least $D+1$ children. Let us denote the set of addition vertices with at least $D+1$ children by $V^+$.
    
    \paragraph{3. Structure of $G(u)$ for $u \in V^+$:}
    Let $u \in V^+$.
    All vertices in $G(u)$ have the same neighborhood $W$ in $V\left(G - G(u)\right)$. Since $|G(u)| \geq t$ it is clear that $w \coloneqq |W| < s$: if $w \geq s$, then $W$ together with any $t$ vertices from $G(u)$ would form a $K_{s,t}$, contradicting $\mathscr{P}_s = t-1$.

    We apply Lemma~\ref{thm:profile_restriction} with $\mathscr{P}^{(1)} = \mathscr{S}(G[W])$, the biclique-sequence of $G[W]$, and obtain a biclique-sequence $\mathscr{P}^{(u)}$ such that $G(u)$ fulfills $\mathscr{P}^{(u)}$ if and only if $G$ fulfills $\mathscr{P}$.
    It is clear, that $G(u)$ is $\mathscr{P}^{(u)}$-extremal.
    By~\eqref{eq:profile_restriction} we have
    \begin{align}
        \mathscr{P}^{(u)}_c = \text{min}_{a \geq 0} \left( \mathscr{P}_{a + c} - \mathscr{S}_a(G[W]) \right).
    \end{align}
    We deduce that $\mathscr{P}^{(u)}$ has start index $s-w$ as for $c < s-w$ either $a \leq w$ so $\mathscr{P}_{a+c} = \infty$ and $\mathscr{S}_a(G[W]) < \infty$ or $a > w$ and $\mathscr{S}_a(G[W]) = - \infty$ (remember that $\mathscr{P}_{a+c}$ has non-negative limit).

    Now we apply Lemma~\ref{thm:lifting} to $G(u)$ with the profile $\mathscr{P}^{(u)}$.
    Each connectivity component of $G(u)$ corresponds to a child $x$ of $u$ in the cotree, which in turn corresponds to a connected cograph $G(x) = G_{x,1} \times G_{x,2}$ with $|G_{x,1}| \leq |G_{x,2}|$.
    By part~(i) of the lemma, $|G_{x,1}| \leq \mathscr{P}^{(u)}_{s-w} = t-1$. We partition the children of $u$ into two sets:
    \begin{itemize}
        \item $X^{(u)}_1$: children $x$ where $|G_{x,1}| \leq s-w-1$. By Lemma~\ref{thm:lifting}~(ii), there is at most one such child for each value $j \in [s-w-1]$, so $|X^{(u)}_1| \leq s-w-1$. The size of $G(x)$ for $x \in X^{(u)}_1$ is unrestricted.
        \item $X^{(u)}_2$: children $x$ where $|G_{x,1}| \geq s-w$. Since also $|G_{x,2}| \geq |G_{x,1}| \geq s-w$, Observation~\ref{obs:component_size_bound} applies to $G(x)$, giving $|G(x)| \leq 2\mathscr{P}^{(u)}_{s-w} = 2(t-1)$.
    \end{itemize}

    \paragraph{4. Reduction to a core graph with one pumped component:}
    For the fixed addition vertex $v \in V^+$ we know that $|X^{v}_2| > D - |X^{v}_2| > D - s = 4t^2$.
    As for $x \in X^{(v)}_2$ we know $|G(x)| < 2t$ and the pigeonhole principle yields $q \in [2t]$ such that there are at least $2t$ vertices $X_{\text{pump}} \subseteq X^{(v)}_2$ such that for any $x \in X'$ we have $|G(x)| = q$.
    
    Let us denote all the children of such addition vertices that correspond to small components except $X_{\text{pump}}$ by
    \begin{align}
        X_{\text{pool}} \coloneqq \bigcup_{u \in V^+} X^{(u)}_2 \setminus X_{\text{pump}}.
    \end{align}

    Now for any $x \in X_{\text{pool}}$ we may choose an arbitrary subset $X_{\text{test}} \subseteq X_{\text{pump}}$ of size $\mathscr{L}(x)$ and replace the vertices $\bigcup_{y \in X_{\text{test}}} \mathscr{L}(y)$ by $q$ times pumping $G(x)$ to obtain a cograph $G'$ of the same vertex count.
    For showing that $G'$ fulfills $\mathscr{P}$ let $u \in V^+$ be the parent of $x$ in the cotree.
    We know that pumping $G(x)$ does not create $\mathscr{P}$-problematic bicliques.
    Let $W$ be the common neighbourhood of $\mathscr{L}(x)$.
    When pumping $G(x)$ new bicliques can only arise in
    \begin{align}
        G[W] \times \left( \sum_{\underset{x' \neq x}{x' \ \text{child of} \ u}} G(x') + m \cdot G(x) \right)
    \end{align}
    However, as $G[W] \times G(x)$ does not contain $\mathscr{P}$-problematic bicliques, the only way such a biclique could arise would be with $W$ and the right hand side acting as the partite sets. However this is not possible as $|G(u)| > D > t$, so $W$ has less than $s$ vertices.    
    We deduce
    \begin{align}
        \| G' \| \leq \| G \|. \label{eq:move_inequality}
    \end{align}

    For $z \in V(T)$ recall the edge contribution $e(\mathscr{L}(z))$ as defined in the Notation section.
    It is clear that for any $x_1, x_2 \in X_{\text{pump}}$ we have $e(\mathscr{L}(x_1)) = e(\mathscr{L}(x_2))$, let us denote it by $e_{\text{pump}}$.
    Hence, we may replace all graphs $G(x)$ for $x \in X_{\text{pump}}$ by one of these components, let us denote it by $G_{\text{pump}}$.
    Moreover, by (\ref{eq:move_inequality}) we have
    \begin{align}
        q \cdot e(\mathscr{L}(x)) \leq |\mathscr{L}(x)| \cdot e_{\text{pump}}.
    \end{align}

    Now, for any $X_{\text{candidates}} \subseteq X_{\text{pool}}$ of size $q$ by Davenport's Theorem \ref{thm:davenport} we find $X_{\text{replace}} \subseteq X_{\text{candidates}}$ such that $q$ divides $\sum_{x \in X_{\text{replace}}} |G(x)|$, let us say by $r$.
    As
    \begin{align}
        q \cdot \sum_{x \in X_{\text{replace}}} e(\mathscr{L}(x))
            \leq e_{\text{pump}} \cdot \sum_{x \in X_{\text{replace}}} |\mathscr{L}(x)|
    \end{align}
    we conclude that
    \begin{align}
        \sum_{x \in X_{\text{replace}}} e(\mathscr{L}(x)) \leq e_{\text{pump}} \cdot r.
    \end{align}
    This shows that we may replace $\bigcup_{x \in X_{\text{replace}}} \mathscr{L}(x)$ by $r$ copies of $G_{\text{pump}}$ without losing edges. By a similar argument as about $G'$ we see that $G_{\text{pump}}$ fulfills $\mathscr{P}$.

    Iterating this procedure, we arrive at a $\mathscr{P}$-extremal cograph $\hat{G}$, where all vertices except $v$ have at most $D$ children, except for at most $q-1$ leftover children from $X_{\text{pool}}$.
    
    Now, let us consider the cograph $G_{\text{core}}$ obtained by deleting $\bigcup_{x \in X_{\text{pump}}} \mathscr{L}(x)$ from $\hat{G}$.
    It has at most $N_{s,t} \coloneqq D^{2(t+s)+1} + q \cdot 2t$ vertices.
    Moreover, it fulfills $\mathscr{P}$.

    \paragraph{5. Cyclic repetition of the core graphs:}
    For any cograph of size at most $N_{s,t}$ in its cotree there is a finite number of children $x$ of addition vertices, so that one can blow up $G(x)$ indefinitely so that the resulting cograph fulfills $\mathscr{P}$.
    This gives a finite number of equation restriction pairs of the form
    \begin{align}
        a_j + \alpha_j \cdot n, \ n \bmod q_j = \gamma_j,
    \end{align}
    for $j \in \Omega$ that correspond to a cograph $G_j$ on at most $N_{s,t}$ vertices together with the child of an addition vertex $u_j \in V(T_{G_j})$ that's corresponding subgraph $G(u_j)$ would be pumped. Observe that $q_j = |\mathscr{L}(u_j)|$ and $\gamma_j = |G_j| \bmod q_j$.
    We conclude that for $n \geq N_{s,t}$
    \begin{align}
        \text{ex}\left( n, \left\{ \mathscr{P}, P_4 - \text{ind} \right\} \right) = \text{max} \left\{ a_j + \alpha_j \cdot n \ | \ j \in \Omega, \ n \bmod q_j = \gamma_j \right\}.
    \end{align}

    As for large $n$ the linear part $\alpha_j \cdot n$ takes over the constant $a_j$ there is $R$ and $\sigma: \ \ZZ \textfractionsolidus R\ZZ \longrightarrow \Omega$ such that for any $n$ large enough
    \begin{align}
        \text{ex}\left( n, \left\{ \mathscr{P}, P_4 - \text{ind} \right\} \right) = a_{\sigma(n \ \text{mod} R)} + \alpha_{\sigma(n \ \text{mod} R)} \cdot n.
    \end{align}

    Moreover, we know that there is $\alpha$ such that for any $j \in \ZZ \textfractionsolidus R \ZZ$ we have $\alpha_{\sigma(j)} = \alpha$. Otherwise, there are distinct $j_1, j_2 \in \ZZ \textfractionsolidus R \ZZ$ such that $\alpha_{\sigma(j_1)} < \alpha_{\sigma(j_2)}$.
    We could add $(R + j_1 - j_2) \bmod R$ singletons to $G_{\sigma(j_2)}$ to obtain $G_{\sigma(j_1)}'$ and pump up the corresponding pumping component appropriatly often.
    This gives $\mathscr{P}$-fulfilling cographs on $n$ vertices with $n \bmod R = j_1$.
    However, when pumping often enough this results in more edges than the original construction using $G_{\sigma(j_1)}$ as in the edge count formula the linear part takes over the constant part for large $n$, a contradiction.

    For notational convenience let us denote $a_j$ for $a_{\sigma(j)}$ and $G_j$ for $G_{\sigma(j)}$ as well as $H_j$ for $H_{\sigma(j)}$.

    \paragraph{6. Upper bound for $\alpha$ and $a_1, \dots,a_{R-1}$.}
    Let us assume for a contradiction that there is $q \in \mathbb{Z}/R\mathbb{Z}$ such that $a_q > 0$.
    Consider a $\mathscr{P}$-extremal cograph $G$ on $n$ vertices with $n \bmod R = q$ that was obtained by pumping up $H_q$ in $G_q$. Now, $(R+1) \cdot G$ still fullfills $\mathscr{P}$ but has more edges than the extremal cograph obtained by pumping $H_q$ to the same vertex count, a contradiction.

    Let us proceed by proving the claimed upper bound of $\alpha$ and $a_0, \dots, a_{R-1}$ by induction on the lexicagraphic order of $(s,t-1)$.
    In case $s=2$ it is clear that $\text{ex} \left( n, \{ \mathscr{P}, P_4 - \text{ind} \} \right) \leq \text{ex} \left( n, \{ K_{2,t}, P_4 - \text{ind} \} \right)$.
    By Theorem \ref{thm:classification_2_t_function}
    \begin{align}
        \text{ex} \left( n, \{ K_{2,t}, P_4 - \text{ind} \} \right) < n \cdot \left( s-1 + \frac{t-1}{2} \right). \label{eq:bound_2_t}
    \end{align}
    This shows that $\alpha \leq s-1 + \frac{t-1}{2}$.
    Moreover, in case there is a pumping component being a connectivity component in its core graph, then $\alpha < s-1 + \frac{t-1}{2}$. Indeed the connectivity component is $\mathscr{P}$-extremal itself and the bound (\ref{eq:bound_2_t}) applies.
    
    Let us assume now that $s \geq 3$ and the statement holds for any biclique-profile with lexicagraphically smaller $(s,t-1)$.
    Let $q \in \mathbb{Z}/R\mathbb{Z}$.

    \paragraph{case 1.}
        $H_q$ is not a connectivity component in $G_q$. Let $W$ be the common neighborhood of $H_q$ in $G_q$ and $w \coloneq |W|$. It is clear that $1 \leq w \leq s-1$.
                
        As we saw in detail in step 3. by Lemma \ref{thm:profile_restriction} we know that there is a biclique-profile $\mathscr{P}'$ with start-index $s-w$ and start value $\mathscr{P}'_{s-w}=t-1$ such that $H_q$ is $\mathscr{P}'$-extremal.
        If $H_q$ would have (internal) edge-density greater than $s-w-1 + \frac{t-1}{2}$, we could consider $m \cdot H_q$ demonstrating that for $\mathscr{P}'$ the linear term of $\text{ex} \left( n, \{ \mathscr{P}', P_4 - \text{ind} \} \right)$ is greater than $s-w-1 + \frac{t-1}{2}$, a contradiction to the induction assumption.
        
        Hence we know that the edge-density of $H_q$ is at most $s - w - 1 + \frac{t-1}{2}$. As any vertex in $H_q$ has exactly $w$ neighbors outside of $H_q$ this immediately implies that $\alpha \leq s -1 + \frac{t-1}{2}$.

        Let again $G$ be an extremal cograph obtained by pumping $H_q$ in $G_q$. Let us assume for a contradiction that $a_q = 0$, i.e. the edge density of $G$ is $\alpha$.
        Let $M$ be the connectivity component of $G$ that contains the pumped $H_q$. It is clear that $M$ is the product of a cograph on at most $s-1$ vertices with another cograph.
        Now, $(R+1) \cdot G$ is a $\mathscr{P}$-extremal cograph by edge count that contains many copies of $M$, a contradiction to the Lifting Theorem \ref{thm:lifting}.

        Hence, we now that the edge density of $G_q$ is less than $\alpha$ and therefore $a_q < 0$.

    \paragraph{case 2.}
        $H_q$ is a connectivity component of $G_q$, i.e. $H_q = M_{1} \times M_{2}$ for two cographs $M_1, M_2$ on $m_1, m_2$ vertices respectively.
        As there are $\mathscr{P}$-extremal cographs where $H_q$ is pumped up by the Lifting Theorem \ref{thm:lifting} we know that $s \leq m_1, m_2 \leq t-1$.
        By Lemma \ref{thm:profile_restriction} know that there is a biclique-profile $\mathscr{P}^{M_1}$ such that $M_1$ is $\mathscr{P}^{M_1}$-extremal. We also know that $M_1$ does not contain $K_{s,t-m_2}$.

        \paragraph{case $s \leq t - m_2$.}
            To show an upper edge bound we may weaken the biclique-profile $\mathscr{P}^{M_1}$ by setting values at indices less than $s$ to infinity.
            In case that the edge density of $M_1$ is at least $s - 1 + \frac{t - m_2 - 1}{2}$, we might clone $M_1$, resulting in a bound $\text{ex}\left( n, \{ \mathscr{P}^{M_1}, P_4 - \text{ind} \} \right) \geq (s - 1 + \frac{t - m_2 - 1}{2}) \cdot n$ for infinitely many $n$, contradicting the induction assumption for $(s,t-m_2-1)$ as $s \geq 2$. We deduce
            \begin{align}
                \text{avgdeg}_{M_1}(M_1) < 2(s-1) + t-m_2-1.
            \end{align}
            Hence, the average degree of vertices $V(M_1)$ inside $H_q$ is upper bounded by
            \begin{align}
                \text{avgdeg}_{H_q}(M_1) < 2(s-1) + t-1.
            \end{align}

        \paragraph{case $t - m_2 \leq s$.}
            To show an upper edge bound we may weaken the biclique-profile $\mathscr{P}'$ by setting values at indices less than $t-m_2$ to infinity. A similar argument using the induction hypothesis for $(t - m_2, s-1)$ yields
            \begin{align}
                \text{avgdeg}_{M_1}(M_1) \leq 2(t - m_2 - 1) + s-1.
            \end{align}
            Hence, the average degree of vertices $V(M_1)$ inside $H_q$ is upper bounded by
            \begin{align}
                \text{avgdeg}_{H_q}(M_1) \leq 2(t - m_2 - 1) + s-1 + m_2 \leq 2(s-1) + t-1,
            \end{align}
            where we used our assumption $t-m_2 \leq s$.
            Similarily to the previous case if $t - m_2 \geq 2$ we also obtain the inequality
            \begin{align}
                \text{avgdeg}_{H_q}(M_1) < 2(s-1) + t-1.
            \end{align}
            
        By symmetry we also see that $\text{avgdeg}_{H_q}(M_2) \leq 2(s-1) + t-1$ with an inequality in case $t-m_1 \geq 2$.
        In case that either $m_1 < t-1$ or $m_2 < t-1$ this shows $\alpha < s-1 + \frac{t-1}{2}$.

        Thus, we may assume that both $m_1=m_2 =t-1$. However, in this case the maximal degree within $M_1$ and $M_2$ is at most $s-1$, resulting in a total edge density of $H_q$ of at most $\frac{s-1 + t-1}{2}$. Again we showed that $\alpha < s-1 + \frac{t-1}{2}$.

        To summarize, either there is a $q$ such that $H_q$ is a connectivity component. In this case, $\alpha < s-1 + \frac{t-1}{2}$. Otherwise, for all $q$ the $H_q$ have common neighbours and we proved that $a_q < 0$.

    \paragraph{7. Upper bound on the extremal function.}

        Let us prove that $\text{ex} \left( n, \{ \mathscr{P}, P_4 - \text{ind} \} \right) < n \cdot \left( s-1 + \frac{t-1}{2} \right)$ for all $n \in \NN$.
        Let us assume for a contradiction that there is $n \in \NN$ such that there is a $\mathscr{P}$-extremal cograph $G$ on $n$ vertices with at least $n \cdot \left( s-1 + \frac{t-1}{2} \right)$ many edges.
        As any of its connectivity components needs to be extremal, we may assume that $G$ is connected.

        Now, we might consider the graphs obtained by pumping up $G$ itself $m$ times.
        By edge count for large $m$ we know that $m \cdot G$ is extremal.
        However, now for $n = m \cdot |G|$ step 6 yields that $\text{ex} \left( n, \{ \mathscr{P}, P_4 - \text{ind} \} \right) = a_{(n \bmod R)} + \alpha \cdot n$ for all $n \in \NN$ for $\alpha < s-1 + \frac{t-1}{2}$ and $a_0, \dots, a_{R-1} \leq 0$. A contradiction.
        This closes the proof of Theorem \ref{thm:pumping_profile}.
\end{proof}

\beginProofOfTheoremClickable{thm:pumping}
    Let us define $\mathscr{P}$ by
    \begin{align}
        \mathscr{P}_j \coloneq \begin{cases}
            \infty \qquad &j < s \\
            t-1 \qquad &s \leq j < t \\
            s-1 \qquad &j \geq t
        \end{cases}
    \end{align}

    Theorem \ref{thm:pumping_profile} yields $R, N_{s,t} \in \NN$ and core graphs $G_0, \ldots, G_{R-1}$ each of size at most $N_{s,t}$ with designated components $H_j \subseteq G_j$ of edge density $\alpha$ such that for any $n$ large enough, there is a $\mathscr{P}$-extremal cograph on $n$ vertices obtained by pumping $H_{(n \bmod R)}$ inside of $G_{(n \bmod R)}$.
    For $r \in \NN_0$ let us define
    \begin{align}
        \tilde{G}(r) \coloneq K_{s-1} \times (r \cdot K_{t}).
    \end{align}
    It is clear that $\tilde{G}(r)$ does not contain $K_{s,t}$ and
    \begin{align}
        n_r &\coloneq | \tilde{G}(r) | = s-1 + r \cdot t \\
        \| G(r) \| &= \binom{s-1}{2} + \left( s-1 + \frac{t-1}{2} \right) \cdot (n_r - s + 1)
    \end{align}
    This proves that $\alpha \geq s-1 + \frac{t-1}{2}$ and by the upper bound given by Theorem \ref{thm:pumping_profile} we conclude
    \begin{align}
        \alpha = s-1 + \frac{t-1}{2}.
    \end{align}

    We know from the proof of Theorem \ref{thm:pumping_profile} that the pumping components $H_q$ can not be connectivity components in $G_q$ as otherwise we would have an inequality here.
    Hence, $H_q$ has a common neighborhood $W$ in $G_q$ of size $w \coloneq |W|$.
    It is clear that $1 \leq w \leq s$.
    Let us assume for a contradiction that $w \leq s-2$.
    Then by Lemma \ref{thm:profile_restriction} there is a biclique-profile $\mathscr{P}'$ with start-index $s - w \geq 2$ and start value $t-1$ such that $m \cdot H_q$ is $\mathscr{P}'$-extremal for any $m \in \NN$.
    However, again Theorem \ref{thm:pumping_profile} yields for the linear coefficient $\alpha'$ of $\text{ex} \left( n, \{ \mathscr{P}', P_{4} - \text{ind} \} \right)$ since $H_q$ is a connectivity component in $m \cdot H_q$ that $\alpha' < s-w-1 + \frac{t-1}{2}$, a contradiction to $\alpha = s-1 + \frac{t-1}{2}$.

    This proves that any vertex in $H_q$ has $s-1$ vertices outside of $H_q$ and $H_q$ $(t-1)$-regular.
    By the equality $\alpha = s-1 + \frac{t-1}{2}$ again Theorem \ref{thm:pumping_profile} yields that the constants $a_0, \dots, a_{R-1} < 0$. Also, as $s \geq 2$ it directly implies for any $n \in \NN$ that
    \begin{align}
        \text{ex} \left( n, \{ K_{s,t}, P_4 - \text{ind} \} \right) < n \cdot \left( s-1 + \frac{t-1}{2} \right).
    \end{align}

    Let $q \in \mathbb{Z}/R\mathbb{Z}$. Then for large $n$ with $n \bmod R = q$ one can obtain an $(s,t)$-extremal cograph $G$ on $n$ vertices by pumping up $H_q$ in $G_q$.
    Let $G_q'$ be the connectivity component of $G_q$ that contains $H_q$. Then in $G$ the connectivity component obtained from pumping $H_q$ in $G_q'$ is itself $(s,t)$-extremal.
    Hence for any $n$ large enough such that $n + |G_q - G_q'| \bmod R = q$ there is an $(s,t)$-extremal connected cograph on $n$ vertices and one could pick at least one residue class such that the respective extremal cographs are connected.
    This closes the proof of Theorem \ref{thm:pumping}.
\end{proof}

\section{Conclusion}

We have shown that the bipartite Turán problem on cographs admits an exact linear solution for large $n$, characterized by a periodic "pumping" of specific dense components.
This contrasts with the general bipartite Turán problem, where bounds are polynomial but not linear. The strong structural constraints of cographs (bounded clique-width, bounded VC-dimension) force the extremal cographs into a rigid algebraic shape.
While it was known that the $\text{ex}\left( n ,\{ K_{s,t}, P_4 - \text{ind} \} \right)$ can be linearily bounded in $n$ by a result of Bonamy et al. \cite{Bonamy2020}, we are the first recognizing the cyclic behaviour and determining the linear coefficient of $\text{ex}\left( n ,\{ K_{s,t}, P_4 - \text{ind} \} \right)$ for large $n$ to be $s-1 + \frac{t-1}{2}$.

As one ingredient to our Pumping Theorem is the Davenport Theorem \ref{thm:davenport}, we could not reduce the core even further: in order to replace the bad structures we need multiple of them to appear in order to fulfill some divisibility condition. On the other end, we develop an algorithm to determine extremal cographs for small $n$. However, even for small parameters $s$ and $t$ the two theoretically and computationally analyzed ranges do not meet. We need other tooling to understand this middle range. As the results for the small and the large range mostly agree, it is natural to assume that also in this middle range pumping like behaviour occurs.

Finally, we want to remark that the exact determination of $\text{ex}\left( n, \left\{ H, P_4 - \text{ind} \right\} \right)$ is still an unsolved and interesting problem for $H$ that is neither a clique nor a biclique.

\section*{Acknowledgements}

The author thanks Maria Axenovich for guiding him into extremal graph theory research and sparking interest in the topic of the induced Turán problem. He is grateful for remarks about related Ramsey-type questions, leading to Corollary \ref{thm:balanced_bicliques}.
The author wants to thank Alexandra Wesolek for encouraging engagement with the problem and for many helpful discussions.
He is grateful to Olaf Parczyk for detailed discussions of the proofs and valuable feedback on earlier drafts of this paper.

\bibliographystyle{plainnat}
\bibliography{references}

\end{document}